\documentclass{article}

\usepackage{amsmath}  \usepackage{amssymb}

\newcommand{\ul}[1]{\underline{#1}}  \newcommand{\lra}{\longrightarrow}

\newcommand{\B}{\mathcal{B}}  \newcommand{\E}{\mathcal{E}}  \newcommand{\K}{\mathcal{K}} 

\newcommand{\U}{\mathcal{U}}

\newcommand{\Z}{\mathbb{Z}}    \newcommand{\R}{\mathbb{R}} 

\newcommand{\RP}{\mathbb{RP}}  \newcommand{\C}{\mathbb{C}}

\renewcommand{\labelenumi}{\arabic{enumi})}

\newcommand{\Supp}{\mbox{Supp}}

\title{On the Baum-Connes conjecture in the real case}  \author{by Paul Baum and Max Karoubi} 

\date{}

\begin{document}

\maketitle

The classical Baum-Connes conjecture (for a given discrete countable group $\Gamma$) states that

the index map \cite{Baum}  $$\mu(\Gamma) : K_j^\Gamma(\ul{E}\Gamma) \lra K_j(C_r^*(\Gamma))$$  is

an isomorphism (where $j=0,1 \mod 2$). 

In this statement, $K_j(C_r^*(\Gamma))$ is the $K$-theory of the reduced $C^*$-algebra
$C_r^*(\Gamma)$ (also denoted $C_r^*(\Gamma ; \C)$ in \cite{Roe3}) and $K_j^\Gamma(\ul{E}\Gamma)$
is the \emph{complex} equivariant Kasparov $K$-homology (with $\Gamma$-compact supports) of the
space $\ul{E}\Gamma$.  This index map may be also defined in the \emph{real} context, by using real
Kasparov theory.  In other words, there is an index map  $$\mu_\R(\Gamma) :
KO_j^\Gamma(\ul{E}\Gamma) \lra K_j(C_r^*(\Gamma ; \R))$$  where $j$ takes its values in $\Z \mod
8$.  We may now ask whether $\mu_\R(\Gamma)$ is also an isomorphism. 

One source of interest in this question (for a given group $\Gamma$) is the result of S. Stolz
(with contributions from J. Rosenberg, P. Gilkey and others): the injectivity of $\mu_\R(\Gamma)$
implies the stable Gromov-Lawson-Rosenberg conjecture \cite{Gromov} about the existence of a
Riemannian metric of positive scalar curvature on 
compact connected spin manifolds with $\Gamma$ as fundamental
group \cite{Schroder}. 

\bf The purpose of this paper is to show that the Baum-Connes conjecture in the real case follows
from the usual (i.e. complex) case. \rm  More precisely, our theorem is the following: 

\bf THEOREM. \it Let $\Gamma$ be a discrete countable group.  If $\mu(\Gamma)$ is an isomorphism
then $\mu_\R(\Gamma)$ is also an isomorphism. \rm 

The proof relies on an interpretation of the index maps $\mu(\Gamma)$ and $\mu_\R(\Gamma)$ as
$K$-theory connecting homomorphisms associated to exact sequences of (real or complex)
$C^*$-algebras \cite{Higson}\cite{Roe2} and also on a general theorem for Banach algebras which
follows directly from a ``descent theorem'' in topological $K$-theory: 

\bf THEOREM \rm \cite{Karoubi}. \it Let $A$ be a Banach algebra over the real numbers and let
$A'=A\otimes_\R\C$ be its complexification.  If $K_i(A')=0$ for all $i \in \Z \mod 2$, then
$K_j(A)=0$ for all $j \in \Z \mod 8$. \rm 

\bigskip 

\bf \large 1. Definition of $\mu(\Gamma)$ and $\mu_\R(\Gamma)$ \normalsize \rm 

\bf 1.1. \rm In this section, we recall the basic definitions of \cite{Baum} and observe that
these definitions extend quite immediately to the real case. 

The universal proper $\Gamma$-space is denoted by $\ul{E}\Gamma$ and $K_j^\Gamma(\ul{E}\Gamma)$
denotes the following colimit  $$\underset{\Delta}{\mbox{colim}} \quad
KK_\Gamma^j(C_0(\Delta),\C)$$  where $\Delta$ runs over all $\Gamma$-compact subspaces of
$\ul{E}\Gamma$ (by definition, a $\Gamma$-subspace $\Delta$ is called $\Gamma$-compact if the
quotient space $\Delta / \Gamma$ is compact).  The composition of the two homomorphisms 
$$KK_\Gamma^j(C_0(\Delta),\C) \lra KK_\Gamma^j(C_0(\Delta) \rtimes \Gamma, C_r^*(\Gamma)) \lra
KK^j(\C, C_r^*(\Gamma))$$  induces (by taking the colimit) the map $\mu$ referred to in the
introduction.  Here the first homomorphism is Kasparov's descent map \cite{Kasparov} and the
second one is induced by the Kasparov product with  $$1 \in KK^0(\C,C_0(\Delta) \rtimes \Gamma) =
K_0(C_0(\Delta) \rtimes \Gamma).$$ 

\bf 1.2 Remark. \rm In this definition of $\mu$, the specific space $\ul{E}\Gamma$ does not play a
particular role.  In other words, if $X$ is \emph{any} proper $\Gamma$-space, we could define in the
same way an ``index map''  $$\mu(X,\Gamma) : K_j^\Gamma(X) \lra K_j(C_r^*(\Gamma)).$$ 

\bf 1.3 Remark. \rm These definitions extend immediately to the real case.  Hence, there is a real
index map  $$\mu_\R(X,\Gamma) : KO_j^\Gamma(X) \lra K_j(C_r^*(\Gamma ; \R)).$$ 

The \emph{real} Baum-Connes conjecture for the group states that
$\mu_\R(\ul{E}\Gamma,\Gamma)=\mu_\R(\Gamma)$ is an isomorphism for all $j \in \Z \mod 8$. 

\bigskip 

\bf \large 2. Index maps and connecting homomorphisms in $K$-theory  \normalsize \rm

\bf 2.1. \rm The strategy for proving our theorem is as follows.  We will describe (in this section)
a $C^*$-algebra whose $K$-theory (real or complex) vanishes precisely when the corresponding
version of the Baum--Connes conjecture is true\footnote{In the analogous context of surgery theory,
the $K$-theory of this $C^*$-algebra would be the `fiber of assembly' or `structure set' term in
the surgery exact sequence.}  In the next section we will apply to this $C^*$-algebra the result
of~\cite{Karoubi}, according to which the $K$-theory of a real $C^*$-algebra vanishes if and only
if the $K$-theory of its complexification vanishes.

\bf 2.2. \rm 
To construct the required $C^*$-algebra we have chosen to use the method of 
\cite{Higson}, \cite{Roe1},
\cite{Roe2}. 
Let $X$ be a locally
compact space\footnote{We assume $X$ to be also second countable in order to get separable Hilbert
spaces.} $X$ and let $\Gamma$ be
a countable discrete group $\Gamma$ acting properly on $X$.
Choose a separable Hilbert
space $H$ with a representation of $C^*$-algebras $\psi : C_0(X) \lra \B(H)$ and a unitary group
representation $\tau : \Gamma \lra \U(H)$ which are compatible in the sense that
 $\psi(\gamma.f) =
\tau(\gamma).\psi(f).\tau(\gamma)^*$, where $\gamma.f$ is the function $x \mapsto
f(\gamma^{-1}x)$.   
(Note that these conditions imply that we have in fact a representation of the crossed product
$C^*$-algebra $C_0(X) \rtimes \Gamma$ in $\B(H)$.)   
It is also required that $H$ be a `large' representation in a certain technical sense; it is
sufficient to take $H=L^2(X;\mu)\otimes\ell^2(\Gamma)\otimes H'$, where $H'$ is an auxiliary
infinite-dimensional Hilbert space and $\mu$ is a Borel measure on $X$ whose support is all of $X$.

Within this setting, we define the support in $X \times X$ of an operator $T$, denoted $\Supp(T)$,
as the complement of the points $(x,y)$ such that there exists a neighborhood $U \times V$ of
$(x,y)$ such that $\psi(f)T\psi(g)=0$, for $f$ supported in $U$ and $g$ supported in $V$. 

\bf 2.3. \rm Following \cite{Higson} and \cite{Roe1}, we define now the $C^*$-algebra
$D_\Gamma^*(X)$ and a closed ideal $C_\Gamma^*(X)$.  Thus there is an exact sequence of
$C^*$-algebras  $$0 \lra C_\Gamma^*(X) \lra D_\Gamma^*(X) \lra D_\Gamma^*(X) / C_\Gamma^*(X) \lra 
0.\qquad\qquad(\E)$$ 

By definition, $D_\Gamma^*(X)$ is the closure of the algebra in $\B(H)$ consisting of all the
(bounded) operators $T$ such that  \begin{enumerate}  \item $T$ is $\Gamma$-invariant, i.e.
$T.\tau(\gamma)=\tau(\gamma)T$ for all $\gamma$ in $\Gamma$.  \item $\Supp(T)$ is
$\Gamma$-compact, i.e. its quotient\footnote{Here $\Gamma$ is acting on $X \times X$ by the
diagonal action.} by $\Gamma$ is compact in $(X \times X) / \Gamma$.  \item For all $f$ in
$C_0(X)$, $T\psi(f)-\psi(f)T$ is a compact operator on $H$.  \end{enumerate} 

The ideal $C^*_\Gamma(X)$ is the closure of the algebra in $\B(H)$ consisting of the (bounded)
operators $T$ which satisfy (1), (2), and a stronger condition  \begin{enumerate} 
\setcounter{enumi}{2}  \renewcommand{\labelenumi}{\arabic{enumi}$'$)}  \item For all $f$ in
$C_0(X)$, $T\psi(f)$ and $\psi(f)T$ are compact operators on $H$. 
\renewcommand{\labelenumi}{\arabic{enumi})}  \end{enumerate} 

\bf 2.4 Example. \rm If $\Gamma$ is a finite group and $X$ is compact, it is well known that the
$K$-theory of the $C^*$-algebra $D_\Gamma^*(X) / C_\Gamma^*(X)$ is the $K$-homology, with a shift
of dimension, of the cross-product algebra $C(X) \rtimes \Gamma$ (this is ``Paschke duality''
\cite{Paschke}).  In the simplest case when $X$ is a point, the exact sequence above is
essentially equivalent to a direct sum of exact sequences of the form  $$0 \lra \K \lra \B(H) \lra
\B(H) / \K \lra 0$$  as many as the number of conjugacy classes in $\Gamma$. 

\bf 2.5 THEOREM \rm \cite{Roe2}. \it For any proper cocompact $\Gamma$-space $X$, there is a
canonical Morita equivalence between the $C^*$-algebra $C_\Gamma^*(X)$ and $C_r^*(\Gamma)$, the
reduced $C^*$-algebra of the group $\Gamma$. \rm 

\bf 2.6 THEOREM \rm \cite{Higson}\cite{Paschke}. \it For any proper $\Gamma$-space $X$, there is a
natural isomorphism  $$K_j^\Gamma(X) := \underset{\Delta}{\mbox{colim}} \quad
KK_\Gamma^j(C_0(\Delta),\C) \overset{\cong}{\lra} K_{j+1}(D_\Gamma^*(X) / C_\Gamma^*(X))$$  where
$\Delta$ runs over all the $\Gamma$-compact subspaces of $X$. \rm 

\bf 2.7 THEOREM \rm \cite{Roe2}. \it For any proper $\Gamma$-space $X$, we have a commutative
diagram \rm  $$\begin{array}{ccc}  K_j^\Gamma(X) & \overset{\mu}{\longrightarrow} &
K_j(C_r^*(\Gamma)) \\  \cong \downarrow & & \cong \downarrow \\ 
K_{j+1}(D_\Gamma^*(X)/C_\Gamma^*(X)) & \overset{\delta}{\longrightarrow} & K_j(C_\Gamma^*(X)) \\ 
\end{array}$$  where $\mu$ is the \rm Baum-Connes \it map and where $\delta$ is the $K$-theory
connecting homomorphism associated to the exact sequence $(\E)$ above. \rm 

\bf 2.8 Remark. \rm It is important to notice that the three theorems above are also true in the
real case (see \cite{Roe3} for a detailed account of this ``real Paschke duality'').  In this
case, $C_r^*(\Gamma)$ has to be replaced by $C_r^*(\Gamma ; \R)$.  The real analogs of the
$C^*$-algebras $D_\Gamma^*(X)$ and $C_\Gamma^*(X)$ shall be denoted $D_\Gamma^*(X;\R)$ and
$C_\Gamma^*(X;\R)$. 

\bf 2.9 COROLLARY. \it The \rm Baum-Connes \it map $\mu : K_j^\Gamma \longrightarrow
K_j(C_r^*(\Gamma))$ is an isomorphism for all $j \in \Z \mod 2$ if and only if the $K$-groups
$K_j(D_\Gamma^*(X)) = 0$ for all $j$.  In the same way, the real \rm Baum-Connes \it map $\mu_\R :
KO_j^\Gamma(X) \longrightarrow K_j(C_\Gamma^*(X;\R))$ is an isomorphism for all $j \in \Z \mod 8$
if and only if the $K$-groups $K_j(D_\Gamma^*(X;\R))=0$ for all $j$. \rm 

\bigskip 

\bf \large 3. Proof of the Baum-Connes conjecture in the real case for a given group $\Gamma$
(assuming its validity for $\Gamma$ in the complex case) \normalsize \rm 

\bf 3.1. \rm As we have shown in the second section, the complex (resp. real) Baum-Connes
conjecture is equivalent to the vanishing of the $K$-groups $K_j(D_\Gamma^*(X))$ (resp.
$K_j(D_\Gamma^*(X;\R))$) for $X=\ul{E}\Gamma$.  If we put $A=D_\Gamma^*(X;\R)$, its
complexification $A'=A \otimes_\R \C$ is isomorphic to $D_\Gamma^*(X)$.  The scheme of the
argument is then the following, where $BC(\Gamma)$ (resp. $BC_\R(\Gamma)$) stands for the
Baum-Connes conjecture (resp. the real Baum-Connes conjecture) for a given discrete group
$\Gamma$:  $$BC(\Gamma) \Longleftrightarrow K_*(A')=0 \Longrightarrow K_*(A)=0 \Longleftrightarrow
BC_\R(\Gamma).$$ 

\bf 3.2. \rm The only point to show is the implication $K_*(A')=0 \Longrightarrow K_*(A)=0$, which
follows from the descent theorem stated in \cite{Karoubi} in the general framework of Banach
algebras.  More precisely, let $A$ be any Banach algebra over the real numbers and $A'$ denote its
complexification $A \otimes_\R \C$.  There is then a cohomology spectral sequence with $E_2^{pq} =
H^p(\RP_2 ; K_{-q}(A'))$ converging to $K_{-q-p}(A) \oplus K_{-q-p+4}(A)$, where $\RP_2$ is the
real projective plane and $H^p$ means usual singular cohomology with local
coefficients\footnote{In fact, there is at most one non zero differential, therefore
$E^3=E^\infty$.}.  The hypothesis $K_*(A')=0$ implies that the $E_2$ term of the spectral sequence
is 0.  Therefore the $E_\infty$ term is also 0.  Since moreover the filtration is finite (because
$\RP_2$ is finite dimensional), $K_*(A)$ must be also 0.

\bigskip \bigskip 

\bf Paul Baum \rm (baum@math.psu.edu): Department of Mathematics, Penn State University,
University Park, PA, 16802, USA. 

\bf Max Karoubi \rm (karoubi@math.jussieu.fr): UFR de Math\'ematiques, Universit\'e Paris 7, 2
place Jussieu, 75251 Paris cedex 05, France. 

\end{document}